\documentclass[11pt]{article}
\usepackage{epsfig}
\newcommand{\comment}[1]{}

\newcommand{\reals}{\mbox{$\mathbb R$}}

\newcommand{\nats}{\mbox{$\mathbb N$}}

\newcommand{\power}{\mbox{$\mathbb P$}}

\newcommand{\dg}{{\mathop{\mathrm{deg}}\nolimits}}

\usepackage{amssymb}

\def\squarebox#1{\hbox to #1{\hfill\vbox to #1{\vfill}}}
\def\qed{\hspace*{\fill}
        \vbox{\hrule\hbox{\vrule\squarebox{.667em}\vrule}\hrule}\smallskip}
\newenvironment{proof}{\begin{trivlist}
  \item[\hspace{\labelsep}{\em\noindent Proof.~}]
  }{\qed\end{trivlist}}

\newtheorem{lemma}{Lemma}[section]
\newtheorem{theorem}[lemma]{Theorem}
\newtheorem{corollary}[lemma]{Corollary}
\newtheorem{proposition}[lemma]{Proposition}
\newtheorem{claim}[lemma]{Claim}
\newtheorem{observation}[lemma]{Observation}
\newtheorem{definition}[lemma]{Definition}

\def\squareforqed{\hbox{\rlap{$\sqcap$}$\sqcup$}}
\def\qed{\ifmmode\squareforqed\else{\unskip\nobreak\hfil
\penalty50\hskip1em\null\nobreak\hfil\squareforqed
\parfillskip=0pt\finalhyphendemerits=0\endgraf}\fi}

\newlength{\tablength}
\newlength{\spacelength}
\newcommand{\tabstar}{\hspace*{\tablength}}
\newcommand{\spacestar}{\hspace*{\spacelength}}
\def\obeytabs{\catcode`\^^I=\active}
{\obeytabs\global\let^^I=\tabstar}
{\obeyspaces\global\let =\spacestar}
\newenvironment{display}{\begingroup\obeylines\obeyspaces\obeytabs}{\endgroup}
\newenvironment{prog}{\begin{display}\parskip0pt\sf}{\end{display}}

\title{On Minkowski sums of simplices}

\author{
{\sl Geir Agnarsson} $^{*}$
\and
{\sl Walter D.~Morris}
\thanks{Department of Mathematical Sciences,
George Mason University,
MS 3F2,
4400 University Drive,
Fairfax, VA -- 22030, USA,
{\tt $\{$geir@math.gmu.edu,wmorris@gmu.edu$\}$}}
}

\comment{
\and
{\sl Bernd Sturmfels..?}
\thanks{ Department of Mathematics,
University of California
Berkeley, CA -- 94720, USA,
{\tt  bernd@math.berkeley.edu}}
}

\date{}

\begin{document}

\maketitle

\begin{abstract}
We investigate the structure of the
Minkowski sum of standard simplices in ${\reals}^r$.
In particular, we investigate the one-dimensional structure,
the vertices, their degrees and the edges in the Minkowski
sum polytope.

\vspace{3 mm}



\vspace{2 mm}

\end{abstract}

\section{Introduction and Definitions}
\label{sec:intro}

Let $[r] = \{1,2,\ldots,r\}$. The {\em standard simplex}
$\Delta_{[r]}$ of dimension $r-1$ is given by
\[
\Delta_{[r]} = \{ (x_1,\ldots,x_r)\in {\reals}^r : x_i\geq 0 \mbox{ for all $i$ },
x_1+\cdots+x_r = 1\}.
\]
Each subset $F\subseteq [r]$ yields a {\em face} $\Delta_F$ of $\Delta_{[r]}$ given by
\[
\Delta_F = \{ (x_1,\ldots,x_r)\in \Delta_{[r]} : x_i = 0 \mbox{ for } i\not\in F\}.
\]
Clearly $\Delta_F$ is itself a simplex embedded in ${\reals}^r$. If ${\cal{F}}$ is
a family of subsets of $[r]$, then we can form the {\em Minkowski sum} of simplices
\[
P_{\cal{F}} = \sum_{F\in{\cal{F}}}\Delta_F =
\left\{\sum_{F\in{\cal{F}}}x_F : x_F \in \Delta_F \mbox{ for each
} F\in{\cal{F}}\right\}.
\]
If $|F| = 2$ for all $F \in {\cal{F}}$, then the polytope
$P_{\cal{F}}$ is called a {\it graphical zonotope}.  Graphical
zonotopes were studied by West et. al.~\cite{Fisher},~\cite{West},
but several questions about them have gone unanswered.  Minkowski
sums of simplices have more recently been studied by Feichtner and
Sturmfels~\cite{EvaMaria-Bernd}, and by
Postnikov~\cite{Postnikov}.  These later papers focus on a the
case when the collection ${\cal{F}}$ is a {\it building set}, i.e.
${\cal{F}}$ contains all singletons, and has the property that,
for any $F_1, F_2 \in {\cal{F}}$, $F_1 \cap F_2 \neq \emptyset$
implies that $F_1 \cup F_2 \in {\cal{F}}$. This property implies
that the polytope $P_{\cal{F}}$ is simple.  Applications of
Minkowski sums of simplices appear in the paper of Morton et.
al.~\cite{Morton}.  Minkowski sums of simplices have also appeared
in the work of Conca~\cite{Conca} and of Herzog and
Hibi~\cite{HerzogHibi}, under the name transversal polymatroids.

\begin{observation}
\label{obs:Mink-dim}
The dimension of the polytope $P_{\cal{F}}$
is given by $\dim(P_{\cal{F}}) = n - c$ where
\[
n = \left|\bigcup_{F\in{\cal{F}}}F\right| \in [r]
\]
and $c$ is the number of connected components of
$\Delta_{\cal{F}}$, the simplicial complex with facets
$\max({\cal{F}})$.
\end{observation}

\begin{proof}
For each $F\in{\cal{F}}$ present $\Delta_F$ by
$\Delta_F = \{ (x_{F;1},\ldots,x_{F;n}) \in \Delta_{[r]} :
x_{F;i} = 0\mbox{ for } i\not\in F\}$.
We then use the equation $\sum_{i\in F}x_{F;i} = 1$
for each $F$ to obtain $x_{F;\max(F)} = 1 - \sum_{i\in F\setminus\{\max(F)\}}x_{F;i}$ and
eliminate $x_{F;\max(F)}$ in the Minkowski sum. By
then counting the free variables, we have the observation.
\end{proof}
From the following more graph theoretic point of view we also can
consider the following: Let $\Delta_1({\cal{F}})$ be the
1-dimensional skeleton of $P_{\cal{F}}$.

\begin{observation}
\label{obs:Mink-dim-forest}
The dimension of the polytope
$P_{\cal{F}}$ is given by $\dim(_{\cal{F}}) = |E(T_{\cal{F}})|$,
the number of edges in a spanning forest of $\Delta_1({\cal{F}})$.
\end{observation}

A {\em face} of $P_{\cal{F}}$ is a subset of $P_{\cal{F}}$ on which a
linear function is maximized.  A vector $c = (c_1,\ldots,c_r)\in {\reals}^r$
defines a partition $C = (C_1,C_2,\ldots,C_s)$ of $[r]$ into nonempty
subsets, so that $c_{i_1} = c_{i_2}$ when $i_1$ and $i_2$ are in
the same part of the partition, and $c_{i_1} < c_{i_2}$ whenever
$i_1 \in C_{\ell_1}, i_2 \in C_{\ell_2}, \ell_1 < \ell_2$.  Then
the points of the face $Q$ that maximizes $c^Tx$ satisfy the
equations
\[
\sum_{i\in C_{\ell}} x_i = |\{F\in {\cal{F}}: F \cap C_{\ell} \neq
\emptyset, F \cap C_m = \emptyset \mbox{ for } m > \ell\}|.
\]
for $\ell = 1,2, \ldots,s$.  The face that maximizes $c^Tx$ is
therefore the Minkowski sum of the simplices in the family
\[
{\cal{F}}^C:=\{F \cap C_{\ell_F} :F\in {\cal{F}}, F \cap
C_{\ell_F} \neq \emptyset,F \cap C_m = \emptyset \mbox{ for } m
> \ell_F\}
\]
The dimension of the face is determined by the number of connected
components of the simplicial complex $\Delta_{{\cal{F}}^C}$.  If
$\Delta_{{\cal{F}}^C}$ is obtained from $\Delta_{\cal{F}}$ by
splitting one of the components of $\Delta_{\cal{F}}$ in two, then
the corresponding face of $P_{\cal{F}}$ is a facet, and the the
coefficients of the vector $c$ corresponding to $C$ can be assumed
to be 0 and 1.  Therefore, all facets of $P_{\cal{F}}$ are of the
form $\sum_{i \in D} x_i = t$ for some subset $D$ of $[r]$ and
integer $t$.  When $\Delta_{{\cal{F}}^C}$ has exactly one
component of size two, say $\{i,j\}$, and otherwise all isolated
elements, then the corresponding face of $P_{\cal{F}}$ is an edge
parallel to $e_i-e_j$.  Vertices of $P_{\cal{F}}$ are points that
maximize linear functions $c^Tx$ in which all components of $c$
are distinct.  If $c_1<c_2<\cdots<c_r$ then component $v_i$ of the
vertex that maximizes $c^Tx$ equals the number of sets $F$ for
which $i$ is the largest element.  In particular, vertices of
$P_{\cal{F}}$ have integer coordinates.

\section{Minkowski sum of a fixed number of simplices}
\label{sec:sum-of-k}

Suppose that ${\cal{F}}$ consists of $k$
subsets $F_1,F_2,\ldots,F_k$ of $[r]$. For each $i \in [r]$,
define $N_{\cal{F}}(i) = \{j \in [k]:i \in F_j\}$.  Let $A$ be a
subset of $[r]$ so that $N_{\cal{F}}(i_1) = N_{\cal{F}}(i_2)$
whenever $i_1$ and $i_2$ are in $A$.  We would like to show how
the combinatorial type of $P_{\cal{F}}$ can be inferred from that
of $P_{\cal{F}'}$, where ${\cal{F}'}$ is obtained from ${\cal{F}}$
by replacing each appearance of $A$ in a set $F$ by the
one-element set $m = \max(A)$.  Afterward, we will restrict
our attention to families in which all of the $N_{\cal{F}}(i)$ are
distinct.

Every point $y \in P_{\cal{F}'}$ corresponds to the simplex
$\Delta(y) := \{z \in {\reals}^r :z_i = y_i, i \notin A,\sum_{i
\in A}z_i = y_m, z_i \ge 0, i \in A\}$ contained in $P_{\cal{F}}$.
Note that $\Delta(y)$ is $(|A|-1)$-dimensional if $y_m > 0$ and a
point otherwise.  Let ${\cal{F}''}$ be the face of ${\cal{F}}$
where $y_m = 0$.  The combinatorial type of $P_{\cal{F}}$ is
therefore that of $\Delta_A \times P_{\cal{F}'}$, with (if
$P_{\cal{F}''}$ is nonempty) the face $\Delta_A \times
P_{\cal{F}''}$ collapsed to a copy of $P_{\cal{F}''}$.  In the
case that $|A|=2$, $P_{\cal{F}}$ is a wedge over $P_{\cal{F}'}$
with foot $P_{\cal{F}''}$.

{\sc Example} Consider the family ${\cal {F}} =
\{\{1,2,3\},\{1,2,4\}\}$ of subsets of $[4]$.  Then
$N_{\cal{F}}(i) = \{1,2\}$ for all $i$ in  $A = \{1,2\}$. The
polytope $P_{\cal{F}}$ is drawn in Figure 1.  The polytope
$P_{\cal{F}'}$ is the two-dimensional cube that is the top face of
the drawing. $P_{\cal{F}''}$ is the vertex $(0,0,1,1)$.
\begin{figure}[tb]
\begin{center}
\leavevmode
\hbox{%
\epsfxsize=4.5in \epsffile{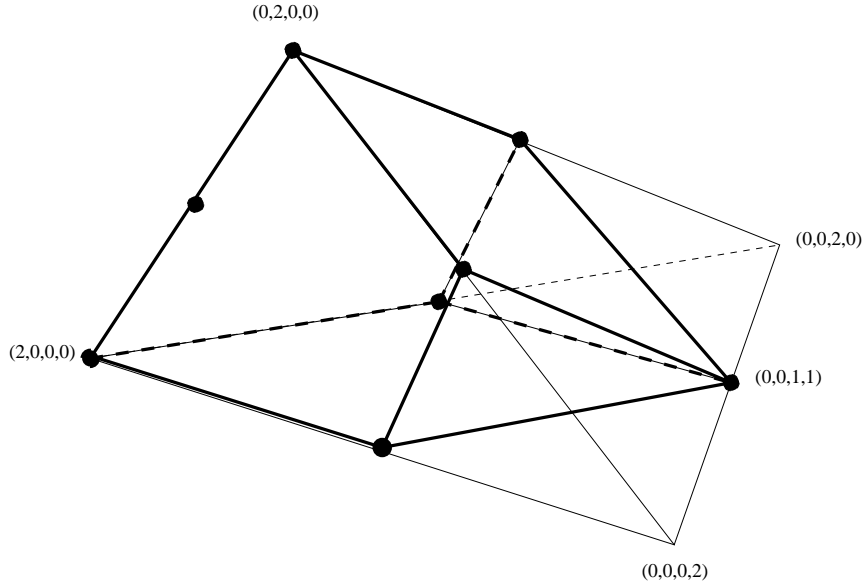}} \caption{A sum of two
triangles}
\end{center}
\end{figure}
\begin{proposition}
\label{prop:wedge}
Every vertex of $P_{\cal{F}}$ is of the form $y + y_m(e_i -
e_m)$, where $y$ is a vertex of $P_{\cal{F}'}$ and $i \in
\{1,2,\ldots,m\}$.  Vertices $y + y_m(e_i - e_m)$ and $y' +
y'_m(e_j -e_m)$ of $P_{\cal{F}}$ are adjacent in $P_{\cal{F}}$ if
\begin{enumerate}
  \item$y = y'$ and $y_m>0$ or
  \item $y$ is adjacent to $y'$ in $P_{\cal{F}'}$ and either $i=j$ or $y_my'_m=0$.
\end{enumerate}
\end{proposition}
Note that no vertex has more than one component
of $A$ nonzero, because the vertices of the simplex $\Delta_A$
have only one nonzero component.

We consider first the case in which ${\cal{F}}$ consists of two
sets, $F$ and $F'$. In the special case where each of the
sets $F\setminus F'$, $F\cap F'$ and $F'\setminus F$ has exactly
one element, say 1,2 and 3 respectively, then $F = \{1,2\}$ and
$F' = \{2,3\}$ and the Minkowski sum $P = \Delta_{F} +
\Delta_{F'}$ is the convex hull of $(1,1,0), (0,2,0), (0,1,1)$ and
$(1,0,1)$ in ${\reals}^3$, which constitutes a two-dimensional
rhombus within the positive octant of the plane $x + y + z = 2$.

We now argue that the generic Minkowski sum of two simplices
roughly has the structure of such a rhombus, if each of
$F\setminus F', F\cap F'$, and $F'\setminus F$ is nonempty.

By assigning the 1st, 2nd and 3d coordinate axis of ${\reals}^3$
to these parts respectively, we can partition the vertices of $P =
\Delta_{F} + \Delta_{F'}$ in the following ``rhombus''-way: A
vertex $e_i + e_j$ of $P_{\cal{F}}$ is of type $A = (1,1,0)$ if
$i\in F\setminus F'$ and $j\in F\cap F'$, of type $B = (0,2,0)$ if
$i = j \in F\cap F'$, of type $C = (0,1,1)$ if $i\in F\cap F'$ and
$j\in F'\setminus F$ and of type $D = (1,0,1)$ if $i\in F\setminus
F'$ and $j\in F'\setminus F$. Note that the rhombus formed by
$A,B,C$ and $D$ in ${\reals}^3$ has edges $AB, BC, CD$ and $DA$.
With this setup we have the following.
\begin{lemma}
\label{lmm:no-diag}
If $A,B,C$ and $D$ are the points in
${\reals}^3$ as here above and $F\setminus F'$, $F\cap F'$ and
$F'\setminus F$ are all nonempty, then there are no $AC$ nor $BD$
type edges of $P = \Delta_{F} + \Delta_{F'}$.
\end{lemma}
\begin{proof}
The original rhombus does not have $AC$ or $BD$ edges.
\end{proof}
By the above Lemma~\ref{lmm:no-diag} we have the following
corollary that describes the structure of a Minkowski sum of two
standard simplices to be roughly that of the rhombus mentioned
above.
\begin{corollary}
\label{cor:rhombus-structure}
If $F,F\subseteq [r]$ then the
edges, or one-dimensional faces, of $P = \Delta_{F} + \Delta_{F'}$
are of the following types:
\begin{enumerate}
  \item Internal $XX$ edges, where both the endvertices
are of type $X \in \{A,B,C,D\}$.
  \item $XY$ edges, with $XY\in \{AB,BC,CD,DA\}$, where one endvertex
is of type $X$ and the other of type $Y$.
\end{enumerate}
\end{corollary}
\begin{theorem}
\label{thm:Mink-of-two}
Let $F,F'\subseteq [r]$ and let $u$ be a
vertex of the polytope $P_{\cal{F}}$.
\begin{enumerate}
  \item If $u$ is of type $A$, $B$ or $C$, then $\dg(u) = |F\cup F'|-1$.
  \item If $u$ is of type $D$, then $\dg(u) = |F|+|F'|-2$.
\end{enumerate}
\end{theorem}
\begin{proof}
If $u$ is of type $B$, say $u = 2e_i$, then $u$ is adjacent to all
$|F \cap F'|-1$ other vertices of type $B$, and all type $A$ and
$C$ vertices of the form $e_i + e_j$, where $j \in (F\setminus F')
\cup (F'\setminus F)$.  If $u$ is of type $A$, say $u = e_i +
e_j$, with $i \in F \setminus F'$ and $j \in F \cap F'$, then $u$
is adjacent to two kinds of type $A$ vertices: $|F \cap F'|-1$
vertices $e_i + e_k$ with $k \in (F \cap F')\setminus \{j\}$ and
$|F \setminus F'|-1$ vertices $e_k + e_j$ with $k \in F \setminus
(F'\cap \{i\})$. Also, $u$ is adjacent to $|F' \setminus F|$ type
$D$ vertices $e_i + e_k$ with $k \in F' \setminus F$, and finally
$u$ is adjacent to the vertex $2e_j$.  If $u$ is of type $D$, say
$u = e_i + e_j$ with $i \in F \setminus F'$ and $j \in F'
\setminus F$, then $u$ is adjacent to $|(F \setminus F') \cup (F
\setminus F')|-2$ vertices of type $D$ obtained by replacing
either $e_i$ or $e_j$ by an $e_k$ for  $k \in (F \setminus F')
\cup (F \setminus F')$, and $u$ is adjacent to $|F \cap F'|$
vertices of each type $A$ and $C$, obtained by replacing $e_i$ or
$e_j$ by an $e_k$ for $k \in F \cap F'$.
\end{proof}
\begin{corollary}
\label{cor:numb-of-edges}
Let $F,F'\subseteq [r]$ and $P = \Delta_{F} + \Delta_{F'}$.
\begin{enumerate}
  \item The total number of vertices of $P$ is
$|F|\cdot |F'| - |F\cap F'|(|F\cap F'| - 1)$.
  \item The total number of one-dimensional
faces (edges) of $P = \Delta_{F} + \Delta_{F'}$ is given by
\[
\frac{1}{2}\left[ |F\setminus F'|\cdot |F'\setminus F| (|F| + |F'|
- 2) + |F\cap F'|(|F\cup F'| - 1) (|F\setminus F'| + |F'\setminus
F| + 1)\right].
\]
\end{enumerate}
\end{corollary}
\begin{proof}
The number of vertices of degree $|F| + |F'| - 2$ in $P$ is
$|F\setminus F'|\cdot |F'\setminus F|$. By
Theorem~\ref{thm:Mink-of-two} the remaining vertices of $P$ all
have degree $|F\cup F'| - 1$. By the Hand-Shaking Theorem the
total number of edges, or one-dimensional faces, is given as
stated.
\end{proof}
Assuming that $F\cup F'= [r]$, then the maximum value of $|F| +
|F'| - 2$ (provided $F\setminus F'$ and $F'\setminus F$ are
nonempty) is $2r-4$, which occurs when $F = [r-1]$ and $F' =
[r]\setminus \{1\}$. Considering the distribution of the two
possible degrees of a Minkowski sum of two simplices $P =
\Delta_{F} + \Delta_{F'}$, we have the following.
\begin{proposition}
\label{prp:average-2} Let $r\in\nats$ be fixed. If $F,F'\subseteq
[r]$ and $P = \Delta_{F} + \Delta_{F'}$ is of dimension $r-1$,
then the average degree $\overline{\dg}(P)$ satisfies
\[
r-1\leq \overline{\dg}(P) < \frac{10}{9}(r-1).
\]
Moreover, the lower bound is attained iff (i) $F\subseteq F'$,
(ii) $F'\subseteq F$ or (iii) $|F\cap F'|=1$. Also,
$\overline{\dg}(P)/(r-1)$ can become arbitrarily close to $10/9$
for large $r$.
\end{proposition}
\begin{proof}
We introduce the variables $x,y$ and $z$ by $x = |F\setminus F'|$,
$y = |F'\setminus F|$ and $z = |F\cap F'|$. Here we have the
boundary condition $x, y\geq 0$ and $x+y+z = r$, and since $P$ is
assumed to have dimension $r-1$ we have $z\geq 1$ or $0\leq x+y
\leq r-1$. By 
Corollary~\ref{cor:numb-of-edges} and the Hand-Shaking Theorem we
obtain that
\begin{eqnarray*}
\overline{\dg}(P) & = &
  2\frac{|E(\Delta_1({\cal{F}}))|}{|V(\Delta_1({\cal{F}}))|} \\
                  & = &
  \frac{|F\setminus F'|\cdot |F'\setminus F|
  (|F| + |F'| - 2) + |F\cap F'|(|F\cup F'| - 1)
  (|F\setminus F'| + |F'\setminus F|
   + 1)}{|F|\cdot |F'| - |F\cap F'|(|F\cap F'| - 1)} \\
                  & = &
  \frac{xy(2r-2-x-y)+(r-1)(r-x-y)(x+y+1)}{(r-y)(r-x)-(r-x-y)(r-x-y-1)}.
\end{eqnarray*}
As a function of $x$ and $y$ we note that $\overline{\dg}(P) =
\overline{\dg}(x,y)$ is symmetric, has the value of $r-1$ on the
boundary of the triangle bounded by $x = 0$, $y = 0$ and $x+y =
r-1$. By 
Theorem~\ref{thm:Mink-of-two} the value
$\overline{\dg}(x,y)$ is strictly larger than $r-1$ inside the
triangle. If the maximum value of $\overline{\dg}(x,y)$ is
$\overline{\dg}_{\max}(r)$, then $(10r - 13)/9 <
\overline{\dg}_{\max}(r) < 10(r-1)/9$, but
$\overline{\dg}_{\max}(r) - (10r - 13)/9$ tends to zero when $r$
tends to infinity.
\end{proof}
{\sc Remark:} In fact, for any $\epsilon > 0$ there is an $r_0$
such that for any $r\geq r_0$ we have
\[
r-1\leq \overline{\dg}(P) < \frac{10r - 13}{9} + \epsilon.
\]

The {\it f-polynomial} $f_P(q)$ of a $d$-dimensional polytope $P$
is $\sum_{i=0}^d f_iq^i$, where $f_i$ is the number of
$i$-dimensional faces of $P$.  Postnikov~\cite{Postnikov} shows
that $f_{P \times Q}(q) = f_P(q)f_Q(q)$ and gives an elegant
formula for $f_{P_{\cal{F}}}(q)$ in the case that ${\cal{F}}$ is a
building set.  If we assume that $A$, ${\cal{F}'}$ and
${\cal{F}''}$ are as in the discussion preceding Proposition
\ref{prop:wedge}, the $f$-polynomial can be decomposed as follows:
\begin{proposition}
$f_{P_{\cal{F}}}(q) = f_{\Delta_A}(q)f_{P_{\cal{F}'}}(q)
- f_{\Delta_A}(q)f_{P_{\cal{F}''}}(q) + f_{P_{\cal{F}''}}(q)$.
\end{proposition}
In the Example, $f_{P_{\cal{F}}}(q) = 7 + 11q +6q^2 + q^3 =
(2+q)(4+4q+q^2) - (2+q)(1) + 1$.

If $P_{\cal{F}}$ is the sum of two simplices $\Delta_{F}$ and
$\Delta_{F'}$, then one can easily check that $P_{\cal{F}} =
\Delta_{F} \times \Delta_{F'}$ when $|F \cap F'|$ is 0 or 1.
This allows us to describe the $f$-polynomials of sums of two
simplices quite easily, using the proposition with
$A = F\cap F'$.
\begin{corollary}
If ${\cal{F}}=\{F,F'\}$, where $F \cap F'
= \{1,2,\ldots,m\}$, then \[f_{P_{\cal{F}}}(q) = f_{\Delta_{F\cap
F'}}(q)f_{\Delta_{(F\cup m)} \times \Delta_{(F'\cup m)}}(q) -
f_{\Delta_{F\cap F'}}(q)f_{\Delta_F \times \Delta_{F'}}(q) +
f_{\Delta_F \times \Delta_{F'}}(q)\]
\end{corollary}
In particular, the number of vertices of
$P_{\cal{F}}$ is
$|F \cap F'|(|F\backslash F'|+1)(|F'\backslash F|+1) -
|F \cap F'||F\backslash F'||F'\backslash F| +
|F\backslash F'||F'\backslash F| = |F \cap F'|(|F\backslash F'| +
|F'\backslash F|+1)+|F\backslash F'||F'\backslash F|$ which is consistent with
Corollary~\ref{cor:numb-of-edges}.

We will now generalize the results that we obtained for the sum of
two simplices to larger sums.
\begin{definition}
\label{def:type-polytope} For $k\in\nats$ let ${\cal{H}}(k)$ be
the family of $k$ subsets of $[2^k-1]$ so that for $i =
1,2,\ldots,2^k-1$, $N_{{\cal{H}}(k)}(i)$ is the $i^{th}$ (in
lexicographic order) nonempty subset of $[k]$. Then $P(k):=
P_{{\cal{H}}(k)}$ is called the $k^{th}$ {\em master polytope}.
\end{definition}

\begin{definition}
Let ${\cal{F}} = (F_1,\ldots,F_k)$ and let $u$ be a point in
$P_{\cal{F}}$. Then $h_{\cal{F}}(u)$ is the point $v$ in $P(k)$
for which, for $i = 1,2,\ldots,2^k-1$, we set
\[
v_i = \left\{
  \begin{array}{ll}
    \sum_{j:N_{\cal{F}}(j) =  N_{{\cal{H}}(k)}(i)} u_j & \mbox{ if there
    is a $j$ with } N_{\cal{F}}(j) =  N_{{\cal{H}}(k)}(i), \\
    0 & \mbox{ otherwise }
  \end{array}
  \right.
\]
\end{definition}

{\sc Remark:} Another way to look at $v = h_{\cal{F}}(u)$ is as
follows: For ${\cal{F}} = (F_1,\ldots,F_k)$ let $u$ be a point in
$P_{\cal{F}}$ for which $u_iu_j > 0$ implies $N_{\cal{F}}(j) \neq
N_{\cal{F}}(i)$. Then let $h_{\cal{F}}(u)$ be the point $v$ in
$P(k)$ where $v_{{\ell}_i} = u_i$ where ${\ell}_i$ is the unique
element in $[2^k-1]$ with $N_{{\cal{H}}(k)}({\ell}_i) =
N_{\cal{F}}(i)$ for each $i\in [r]$.
\begin{theorem}
\label{thm:vert-gen} For ${\cal{F}} = (F_1,\ldots,F_k)$ the point
$u \in P_{\cal{F}}$ is a vertex of $P_{\cal{F}}$ if, and only if,
the following conditions are met.
\begin{enumerate}
  \item Each instance of
$u_{i_{\alpha}}u_{i_{\alpha}}>0$, $N_{\cal{F}}(i_{\alpha}) = N_{\cal{F}}(i_{\beta})$
implies that $i_{\alpha} = i_{\beta}$.
  \item $h_{\cal{F}}(u)$ is a vertex of the polytope $P(k)$.
\end{enumerate}
\end{theorem}

\begin{proof}
(Theorem~\ref{thm:vert-gen} Sketch) For a point $u = e_{i_1} +
\cdots + e_{i_k}$ of $P_{\cal{F}}$ we first note that if
$N_{\cal{F}}(i_{\alpha}) = N_{\cal{F}}(i_{\beta})$ and
$i_{\alpha}\neq i_{\beta}$, then $u = (v + w)/2$ where $v$ and $w$
are the points of $P_{\cal{F}}$ obtained from $u$ on one hand by
replacing $i_{\alpha}$ by $i_{\beta}$ to get $v$ and on the other
hand by replacing $i_{\beta}$ by $i_{\alpha}$ to get $w$. Hence,
the first condition is necessary.

Assume that $u$ satisfies the first condition and that
$h_{\cal{F}}(u)$ is an extreme point of $P(k)$. Since there is a
supporting hyperplane in ${\reals}^{2^k-1}$ containing
$h_{\cal{F}}(u)$ there is a corresponding supporting hyperplane in
${\reals}^n$ containing $u$, showing that $u$ is a vertex of
$P_{\cal{F}}$.

Assume finally that $u$ satisfies the first condition and that
$h_{\cal{F}}(u)$ is not an extreme point of $P(k)$. In this case
$h_{\cal{F}}(u)$ is a proper convex combination of extreme points
of $P(k)$. Since the first condition is satisfied, there are
corresponding points of $P_{\cal{F}}$, such that $u$ is a proper
(in fact the same!) convex combination of these. This completes
the proof.
\end{proof}
For ${\cal{F}} = (F_1,\ldots, F_k)$ let $A_1,\ldots, A_h$
be the vertices of the polytope $P(k)$.
Similar to the case when $k = 2$
we have the following.
\begin{theorem}
\label{thm:gen-edge-structure} If ${\cal{F}} = (F_1,\ldots, F_k)$,
then the edges, or one-dimensional faces, of $P_{\cal{F}}$ are of
the following types:
\begin{enumerate}
  \item Internal $A_iA_i$ type edges, where both the endvertices
are of type $A_i$ for some $i\in\{1,\ldots,m\}$.
  \item $A_iA_j$ type edges, where $A_iA_j$ is an edge of
the master polytope $P(k)$.
\end{enumerate}
\end{theorem}
\begin{proof}
(Sketch) Similarly to the proof of Lemma~\ref{lmm:no-diag}
(although with a bit more elaborate indexing scheme) one can show
that there is a supporting hyperplane in ${\reals}^n$ of $P$
containing the vertex of type $A_i$ and the vertex of type $A_j$
if, and only if, there is a corresponding supporting hyperplane in
${\reals}^{2^k - 1}$ of $P(k)$ containing the vertices $A_i$ and
$A_j$.
\end{proof}
Theorems~\ref{thm:vert-gen} and~\ref{thm:gen-edge-structure} both
reduce the structure of $P_{\cal{F}} \subseteq {\reals}^n$ to
considerations of the master polytope $P(k) \subseteq {\reals}^{2^k-1}$.

\section{Function Representation of Integer Points of $P_{\cal{F}}$}
\label{sec:Newton}

As in the previous section, we assume that ${\cal{F}} =
(F_1,\ldots,F_k)$, an ordered collection of $k$ subsets of $[r]$.
A function  $f : [k] \rightarrow [r]$ that satisfies $f(i)\in F_i$
for each $i$ will be called a  {\em rep-function}.  For a
rep-function $f$ we define $u(f):=e_{f(1)} + \cdots + e_{f(k)}$.
\begin{claim}
\label{clm:min-max} For functions $f, g : [k]\rightarrow [m]$ we
have
\begin{enumerate}
  \item $u(f) + u(g) = u(\min\{f,g\}) + u(\max\{f,g\})$.
  \item If $f \neq g$, then $u(f) \neq u(\min\{f,g\})$.
\end{enumerate}
\end{claim}
In the case $u(f) = u(g)$, we obtain by Claim~\ref{clm:min-max}
that $u(f) = u(g) = ( u(\min\{f,g\}) + u(\max\{f,g\}))/2$.  Hence,
if an integer point $u \in P_{\cal{F}}$ can be represented by two
distinct functions $f$ and $g$, then it is not a vertex of the
type polytope $P(k)$. The interesting part is the converse.

\begin{lemma}
\label{lmm:integer-pt} If $v$ is an integer point in $P_{\cal F}$
that is not a vertex of $P_{\cal F}$, and an edge of the smallest
face containing $v$ is parallel to $e_{i_1} -e_{i_2}$, then
$P_{\cal F}$ contains the points $v + e_{i_1} - e_{i_2}$ and $v -
e_{i_1} + e_{i_2}$.
\end{lemma}
\begin{proof} First note that $v_{i_1} \neq 0$ and $v_{i_2} \neq
0$, because otherwise all points on the smallest face containing
$v$ would satisfy $x_{i_1}=0$ or $x_{i_2}= 0$, contradicting the
assumption that there is an edge of this face parallel to $e_{i_1}
-e_{i_2}$.  If $v$ is on a facet of $P_{\cal F}$ given by $\sum_{i
\in T} x_i = t$ for some $T \subset[r]$ and integer $t$, then this
equation is satisfied by all points in the smallest face
containing $v$. That means that $i_1$ and $i_2$ are either both in
or both outside of $T$. Thus $v + e_{i_1} - e_{i_2}$ and $v -
e_{i_1} + e_{i_2}$ will satisfy any equations that $v$ satisfies.
Furthermore, any inequality $x_i \ge 0$ or $\sum_{i \in T} x_i \le
t$ that $v$ satisfies strictly will also be satisfied by $v +
e_{i_1} - e_{i_2}$ and $v - e_{i_1} + e_{i_2}$, because only one
component is increased by 1 and one component is decreased by 1.
\end{proof}

\begin{lemma}
\label{lmm:altpath} If $f$ and $g$ are rep-functions and $u(g) =
u(f) + te_{i_1} - te_{i_2}$ for $i \neq j$ in $[r]$, then there
exist rep-functions $f_1,f_2,\ldots f_{t-1}$ so that $u(f) +
le_{i_1} - le_{i_2} = u(f_l)$ for $l = 1,2,\ldots,t-1$.
\end{lemma}
\begin{proof} Define $G_{\cal{F}}$ to be the bipartite graph with
vertex set $\{w_j:j\in [k]\} \cup \{v_t:i\in [r]\}$ and edges
$\{(w_j,v_i)\}$ for all $(i,j)$ with $i \in F_j$.  For any
rep-function $h$, let $M_h$ be the set of edges $(w_j,v_i)$ for
which $h(j) = i$. For every $i \in [r]\{i_1,i_2\}$, the number of
edges of $M_g$ meeting $v_i$ equals the number of edges of $M_f$
meeting $v_i$.  For every $j \in [k]$, $w_j$ is met by exactly one
edge from each of $M_f$ and $M_g$.  On the other hand, $v_{i_1}$
is adjacent to $t$ more edges of $M_g$ than $M_f$, and $v_{i_2}$
is adjacent to $t$ more edges of $M_f$ than $M_g$. There therefore
exists a path $P$ from $v_{i_2}$ to $v_{i_1}$ that alternates
between edges of $M_f$ and $M_g$. Let $M^1$ be the set of edges
obtained from $M_f$ by replacing the edges of $M_f$ in the path by
the edges of $M_g$ in the path.  Then, for $j = 1,2,\ldots,k$,
define $f_1(j)=i$, where $(w_j,v_i)$ is an edge of $M^1$.  Then
$u(f_1)=u(f)+e_{i_1}+e_{i_2}$.  We can continue this way to get
$u(f_2),\ldots,u(f_{t-1})$.
\end{proof}
\begin{proposition}
\label{prop:B=integer} Every integer point $v$ in $P_{\cal F}$ is
$u(f)$ for some rep-function $f$.
\end{proposition}
\begin{proof}
The proof is by induction on the dimension of the smallest face
containing $v$.  From the first section, we know that the
statement is true if true if $v$ is a vertex. Suppose $v$ is not a
vertex. Suppose that there is an edge of the smallest face
containing $v$ that is parallel to $e_{i_1} -e_{i_2}$.  Then
lemma~\ref{lmm:integer-pt} allows us to build a segment parallel
to $e_{i_1} -e_{i_2}$, containing $v$ in its interior, and with
endpoints on faces of $P_{\cal F}$ that are of lower dimension
than the one containing $v$.  By induction, the endpoints of the
interval are $u(f)$ and $u(g)$ for some rep-functions $f$ and $g$.
Lemma \ref{lmm:altpath} then gives us a rep-function for $v$.
\end{proof}
\begin{theorem} An integer point $v$ in $P_{\cal{F}}$ is a vertex of
$P_{\cal{F}}$ if and only if there is a unique rep-function $f$ so
that $u(f) = v$. \end{theorem}
\begin{proof} Let $v$ be an integer point in $P_{\cal{F}}$ that is
not a vertex of $P_{\cal{F}}$.  By Lemma~\ref{lmm:integer-pt}
there are $i_1$ and $i_2$ in $[r]$ so that $P_{\cal F}$ contains
the points $v - e_{i_1} + e_{i_2}$ and $v - e_{i_1} + e_{i_2}$.
Let$f$ and $g$ be the rep-functions guaranteed by
Proposition~\ref{prop:B=integer} for $v - e_{i_1} + e_{i_2}$ and
$v - e_{i_1} + e_{i_2}$, respectively.  Let $G_{\cal{F}}, M_f$ and
$M_g$ be as in the proof of Lemma~\ref{lmm:altpath}.  Then There
are two edges of $M_f$ adjacent to $v_{i_2}$ that are not in
$M_g$.  Therefore we can use these edges as initial edges in two
different paths from $v_{i_2}$ to $v_{i_1}$ that alternate between
edges of $M_f$ and $M_g$.  Swapping edges of $M_f$ for edges of
$M_g$ along each of these alternating paths leads to two different
rep-functions for $v$.
\end{proof}

The number of rep-functions for a given ${\cal{F}}$ is easy to
count, it is $\Pi_{F \in {\cal{F}}}|F|$.  By listing the
rep-functions and the corresponding integer points $u(f)$, and
striking out the $u(f)$ that appear more than once, one can list
the vertices of $P_{\cal{F}}$.  This was done by Bernd
Sturmfels~\cite{Bernd-email} for the polytopes $P(k)$, $k=3,4,5$.
He found that $P(3)$ had 41 vertices, $P(4)$ had 1015 vertices,
and $P(5)$ had 59072 vertices.

\section{Max-degree as function of parameters alone}
\label{sec:r-alone}

In this section we determine the function $d : \nats \rightarrow \nats$
given by
\[
d(r) = \max_{\cal{F}}\left\{\dg_{\max}(P_{\cal{F}})\right\},
\]
where the maximum is taken over all multi-subsets $(F_1,\ldots,F_k)$
of $\power([r])$, where $k\in\nats$ can be any integer but $r$ is fixed.
Moreover, for each fixed $k\in\nats$ we determined the function
$d_k : \nats \rightarrow \nats$ defined by
\[
d_k(r) = \max_{|{\cal{F}}| \leq k}\left\{\dg_{\max}(P_{\cal{F}})\right\},
\]
where the maximum is here taken over all multi-subsets
$(F_1,\ldots,F_k)$ of $\power([r])$ where both $k$ and $r$ are fixed.
Clearly $d(r) = \max_{k\in\nats}\{d_k(r)\}$.

We start with the following lower bound for $d_k(r)$
and $d(r)$.
\begin{lemma}
\label{lmm:d-lower}
For $k,r\in\nats$ we have $d_k(r) \geq k(r-k)$,
and therefore $d(r)\geq \lfloor r^2/4\rfloor$.
\end{lemma}
\begin{proof}
Let $k\in [r]$ and let for each $i\in [k]$ let $F_i =
\{i,k+1,k+2,\ldots,r\}$. Then the vertex $v = e_1+e_2+ \cdots e_k$
is adjacent to each of the vertices $v+(e_{i_1}-e_{i_2})$, for $1
\le i_2 \le{i_1}$ and $k+1 \le i_1 \le r$. Therefore $d_k(r)\geq
k(r-k)$, so we have in particular that $d(r)\geq \lfloor
r/2\rfloor\lceil r/2\rceil = \lfloor r^2/4\rfloor$.
\end{proof}
Another polytope that has vertices of degree  $\lfloor
r^2/4\rfloor$ is the graphical zonotope for the complete bipartite
graph with $\lfloor r/2\rfloor$ vertices on one side of the
bipartition and $\lceil r/2\rceil$ vertices on the other side.
West \cite{West} proved that the graphical zonotope for the
complete bipartite graph has vertices of degree $\ell$ for all
$r-1 \le \ell \le \lfloor r^2/4\rfloor$.  On the other hand, every
vertex of the polytope of lemma \ref{lmm:d-lower} other than $v$
has degree $r-1$.

For a fixed vertex $u$, each edge of $P$ incident to $u$ can be
identified with a multiple of a difference $e_i - e_j$ of some
pair of unit vectors, where $i,j\in[r]$ are distinct. Since the
collection $\{ \alpha(e_i - e_j) : \alpha\in\nats\}$ is a set of
parallel vectors, at most one multiple of $e_i - e_j$ can possibly
correspond to an edge incident to $u$. From this alone we see that
the maximum number of edges incident to $u$ is at most ${r\choose
2}$. However, more can be said:

For a vertex $u$ of $P$, let $\vec{G}(u)$ be the directed
graph with the vertexset $V(\vec{G}(u)) = [r]$ where
a directed edge $(i,j)$ is present iff
$u + \alpha(e_i - e_j)$ is a neighbor of $u$ in $P$
for some $\alpha\in\nats$.
\begin{proposition}
\label{prp:deg=digraph}
For $r\in\nats$ and ${\cal{F}} = (F_1,\ldots,F_k)\subseteq \power([r])$,
the digraph $\vec{G}(u)$ is acyclic and its underlying graph $G(u)$
is simple and triangle-free.
\end{proposition}
\begin{proof}
Assume there is a cycle $(i_1,i_2,\ldots,i_h)$ in $\vec{G}(u)$.
Then $u, v_1,\ldots v_h$
are all vertices of $P$, where
$v_{\ell} = u + \alpha_{\ell}(e_{i_{\ell}} - e_{i_{\ell + 1}})$
(here we compute cyclically, so $e_{i_{h+1}} = e_{i_1}$).
This is however impossible since
\[
\sum_{\ell = 1}^h \frac{1}{\alpha_{\ell}}\left(v_{\ell} - u\right) = 0,
\]
which means that there is no hyperplane containing $u$ alone and having
all the $v_{\ell}$'s strictly on one side of it. In particular for $h = 2$,
there are no directed 2-cycles and hence the underlying graph $G(u)$ is simple.
Also for $h = 3$, there are no directed triangles in $\vec{G}(u)$ either.

Assume now that $G(u)$ has a triangle, which then does not correspond
to a directed triangle in $\vec{G}(u)$, say $v = u + \alpha(e_i - e_j)$,
$v' = u + \beta(e_j - e_l)$ and $v'' = u + \gamma(e_i - e_l)$. In this
case we have
\[
v'' - u = \frac{\gamma}{\alpha}\left(v - u\right) +
\frac{\gamma}{\beta}\left(v' - u\right),
\]
which means that the vector $v'' - u$ is in the cone spanned
by $v - u$ and $v' - u$. This contradicts the fact that $uv''$
is an edge of $P$. Hence, the underlying graph $G(u)$ of $\vec{G}(u)$
has no triangles.
\end{proof}
\begin{theorem}
\label{thm:d-upper}
For $r\in\nats$ we have $d(r)\leq \lfloor r^2/4\rfloor$.
\end{theorem}
\begin{proof}
The maximum degree of a vertex $u$ of $P$ is by
Proposition~\ref{prp:deg=digraph} the maximum number
of edges the simple triangle free graph $G(u)$ can
have. By a theorem by Mantel~\cite{Mantel} (as a special
case of Tur\'{a}n's Theorem~\cite{Turan}), the maximum number
of edges of a simple triangle-free graph on $r$ vertices
is $\lfloor r^2/4\rfloor$, hence the theorem.
\end{proof}
By Lemma~\ref{lmm:d-lower} and Theorem~\ref{thm:d-upper} we
have the following corollary.
\begin{corollary}
For $r\in\nats$ we have $d(r) = \lfloor r^2/4\rfloor$.
\end{corollary}
We now turn our attention to the computation of $d_k(r)$.
Note that the Minkowski sum $P_{\cal{F}}$ provided
in the proof of Lemma~\ref{lmm:d-lower} that attains the
overall maximum degree $d(r)$ has
$k = |{\cal{F}}| = \lfloor r/2\rfloor$. Therefore when computing
$d_k(r)$ we can assume $1\leq k\leq r/2$.

First we need a variation of the theorem by Mantel~\cite{Mantel}:
Let $G$ be a simple graph on $n$ vertices and let $1\leq k\leq n/2$.

Call $G$ a {\em $k$-triangle-free graph}, or
a {\em $k$-tr} for short, if $G$ is triangle free and $G$
has a vertex cover of cardinality at most $k$.
\begin{theorem}
\label{thm:var-Mantel}
Let $n\in\nats$ and $1\leq k\leq n/2$. If $e_k(n)$
is the maximum number of edges of a $k$-tr graph $G$,
then $e_k(n) = k(n-k)$.
Moreover, if $G$ is a $k$-tr graph on $n$ vertices
with $e_k(n)$ edges, then $G$ is a complete bipartite with
parts of cardinalities $k$ and $n-k$.
\end{theorem}
\begin{proof}
For $n \in\{1,2\}$ the theorem is trivial. We proceed by
induction and assume we have a $k$-tr graph on $n>2$
vertices with the maximum number $e_k(n)$ of edges.
Let $uv\in E(G)$ be an edge and since either $u$ or $v$
is in the vertex cover $U$ of size $k$, we assume it to
be $u$. Since $G$ is triangle-free the set of neighbors
$N(u)$ and $N(v)$ are disjoint. Let $G' = G - \{u,v\}$
be the simple graph obtained from $G$ by removing the
vertices $u$ and $v$ from $G$. By the disjointness of
$N(u)$ and $N(v)$ we have
$|E(G)| = |E(G')| + d(u) + d(v) - 1$.

Assume first that $v\in U$. In this case $G'$
is a $(k-2)$-tr graph on $n-2$ vertices and hence
by induction hypothesis we have
$|E(G)| = |E(G')| + d(u) + d(v) - 1 \leq
(k-2)[(n-2) - (k-2)] + n - 1 < k(n-k)$.

Now assume that $v\not\in U$. In this case $G'$
is a $(k-1)$-tr graph on $n-2$ vertices and hence
by induction hypothesis we have
$|E(G)| = |E(G')| + d(u) + d(v) - 1 \leq
(k-1)[(n-2) - (k-1)] + n - 1 = k(n-k)$.
Also by inducting hypothesis, $|E(G)| = k(n-k)$ can hold
iff $G'$ is a complete bipartite graph with parts
of cardinalities $k-1$ and $n-k-1$, and $d(u) + d(v) = n$
(i.e.~$N(u)\cup N(v) = V(G)$). This means that $|E(G)| = k(n-k)$
can hold iff $N(v) = U$ and $N(v) = V(G)\setminus U$, that is,
$G$ is a complete bipartite graph with parts of sizes $k$
and $n-k$. This completes the proof.
\end{proof}
From Theorem~\ref{thm:var-Mantel} we obtain the following
corollary.
\begin{corollary}
For $r\in\nats$ and $k\in \{1,\ldots,\lfloor r/2\rfloor\}$,
we have $d_k(r) = k(n-k)$.
\end{corollary}
\begin{proof}
Consider a point $u = e_{i_1} + \cdots + e_{i_k}$ of $P_{\cal{F}}$
(note that some indices might coincide). As noted before, a neighbor
$v$ of $u$ in $P$ must have the form $v = u + \alpha(e_i - e_j)$
for some $\alpha\in\nats$, and $i\in [r]$ and $j\in \{i_1,\ldots, i_k\}$.
Since each directed edge $(i,j)\in V(\vec{G}(u))$ has its head
in $\{i_1,\ldots, i_k\}$, of cardinality at most $k$, the underlying
graph $G(u)$ has a vertex cover of size at most $k$. Therefore
$G(u)$ is a $k$-tr graph and hence by Theorem~\ref{thm:var-Mantel}
at most $k(r-k)$ edges.

In the proof of Lemma~\ref{lmm:d-lower} an example of
$P_{\cal{F}}$ with $|{\cal{F}}| \leq k$ and a vertex
of degree $k(n-k)$ was given. This completes the argument.
\end{proof}

\section{Minkowski sum of three simplices}
\label{sec:sum-of-3}

In this section we will investigate the polytope $P(3)$. Let
${\cal{H}}:={\cal{H}}(3)=(\{1,2,4,5\},\{1,2,3,6\},\{1,3,4,7\})$.
Henceforth we will drop the $(3)$.  Then
$N_{\cal{H}}(1)=\{1,2,3\}, N_{\cal{H}}(2)=\{1,2\},
N_{\cal{H}}(3)=\{2,3\}, N_{\cal{H}}(4)=\{1,3\},
N_{\cal{H}}(5)=\{1\}, N_{\cal{H}}(6)=\{2\}, N_{\cal{H}}(7)=\{3\}$,
so all of the nonempty subsets of $[3]$ are represented.  The case
of $k = |{\cal{F}}| = 3$ is the first interesting case for the
mere reason that the polytope $P(3)$ does not have $2^{k(k-1)} =
64$ vertices, as was the case for $k = 2$, where the rhombus
$P(2)$ had precisely $2^{k(k-1)} = 4$ vertices.

{\sc Example:} the point $A = (0,1,1,1,0,0,0)$ in $P(3)$ is not a
vertex, because $A = (B + C + D)/3$, where $B = (0,2,1,0,0,0,0)$,
$C = (0,0,2,1,0,0,0)$ and $D = (0,1,0,2,0,0,0)$ and all the points
$B,C$ and $D$ are points in the polytope $P(3)$.

\begin{lemma}
\label{lmm:type-pltp-3} The polytope $P(3)$ has 41 vertices in
${\reals}^7$ given by the column vectors (without the last entry)
in the following $7\times 10$, $7\times 21$ and $7\times 10$
matrices. The last entry in each column is the degree of the
vertex.
\[
\begin{array}{ccccc ccc cc}
3 & 1 & 1 & 0 & 0  &  1 & 0 & 0   & 0 & 0 \\
0 & 0 & 2 & 2 & 1  &  0 & 1 & 2   & 0 & 0 \\
0 & 2 & 0 & 1 & 2  &  0 & 0 & 0   & 2 & 1 \\
0 & 0 & 0 & 0 & 0  &  2 & 2 & 1   & 1 & 2 \\
0 & 0 & 0 & 0 & 0  &  0 & 0 & 0   & 0 & 0 \\
0 & 0 & 0 & 0 & 0  &  0 & 0 & 0   & 0 & 0 \\
0 & 0 & 0 & 0 & 0  &  0 & 0 & 0   & 0 & 0 \\
\hline
6 & 6 & 6 & 6 & 6  &  6 & 6 & 6   & 6 & 6
\end{array}
\]
\[
\begin{array}{ccccccc ccccccc ccccccc} 
2 & 1 & 1 & 0 & 0 & 0 & 0    & 2 & 1 & 1 & 0 & 0 & 0 & 0   & 2 & 1 & 1 & 0 & 0 & 0 & 0 \\
0 & 0 & 1 & 1 & 1 & 0 & 0    & 0 & 0 & 1 & 1 & 1 & 0 & 0   & 0 & 0 & 0 & 1 & 1 & 0 & 2 \\
0 & 0 & 0 & 1 & 0 & 1 & 2    & 0 & 1 & 0 & 1 & 0 & 1 & 0   & 0 & 0 & 1 & 1 & 0 & 1 & 0 \\
0 & 1 & 0 & 0 & 1 & 1 & 0    & 0 & 0 & 0 & 0 & 1 & 1 & 2   & 0 & 1 & 0 & 0 & 1 & 1 & 0 \\
1 & 1 & 1 & 1 & 1 & 1 & 1    & 0 & 0 & 0 & 0 & 0 & 0 & 0   & 0 & 0 & 0 & 0 & 0 & 0 & 0 \\
0 & 0 & 0 & 0 & 0 & 0 & 0    & 1 & 1 & 1 & 1 & 1 & 1 & 1   & 0 & 0 & 0 & 0 & 0 & 0 & 0 \\
0 & 0 & 0 & 0 & 0 & 0 & 0    & 0 & 0 & 0 & 0 & 0 & 0 & 0   & 1 & 1 & 1 & 1 & 1 & 1 & 1 \\
\hline
6 & 6 & 6 & 6 & 8 & 6 & 8    & 6 & 6 & 6 & 8 & 6 & 6 & 8   & 6 & 6 & 6 & 6 & 6 & 8 & 8
\end{array}
\]
\[
\begin{array}{c ccc ccc ccc}
1 & 0 & 0  & 1 & 0 & 0  & 1 & 0 & 0   & 0 \\
0 & 0 & 0  & 0 & 1 & 0  & 0 & 1 & 0   & 0 \\
0 & 1 & 0  & 0 & 0 & 0  & 0 & 0 & 1   & 0 \\
0 & 0 & 1  & 0 & 0 & 1  & 0 & 0 & 0   & 0 \\
1 & 1 & 1  & 0 & 0 & 0  & 1 & 1 & 1   & 1 \\
1 & 1 & 1  & 1 & 1 & 1  & 0 & 0 & 0   & 1 \\
0 & 0 & 0  & 1 & 1 & 1  & 1 & 1 & 1   & 1 \\
\hline
7 & 8 & 8  & 7 & 8 & 8  & 7 & 8 & 8   & 9
\end{array}
\]
\end{lemma}
These computations were verified using the computer program POLYMAKE \cite{Gawrilow}.
Using POLYMAKE, we determined that the polytope $P(4)$ had vertices of all degrees in
the set $\{14,15,\ldots,28\}$ except for $\{16,23,26,27\}$.

\subsection*{Acknowledgments}

The authors would like to thank James F.~Lawrence
for helpful discussions regarding the theory of polytopes
in general. Last but not least, sincere thanks to Bernd Sturmfels for
introducing this problem to us and for his keen encouragement.

\flushright{\today}
\end{document}